\documentclass[a4paper,10pt,reqno]{amsart}
\usepackage[a4paper,tmargin=40mm,bmargin=35mm,hmargin=30mm]{geometry}
\usepackage[T1]{fontenc}
\usepackage{lmodern}
\usepackage{amsmath}
\usepackage{amssymb}
\usepackage{amsthm}
\usepackage{stmaryrd}
\usepackage{tikz}
\usetikzlibrary{cd}
\usepackage{booktabs}
\usepackage{multirow}
\usepackage{multicol}
\usepackage{enumitem}
\usepackage{here}
\usepackage{aliascnt}
\usepackage{hyperref}
\usepackage[noabbrev]{cleveref}

\newcommand{\NewTheorem}[2]{
	\newaliascnt{#1}{TheoremEnvironment}
	\newtheorem{#1}[#1]{#1}
	\aliascntresetthe{#1}
	\crefname{#1}{#1}{#2}
	\Crefname{#1}{#1}{#2}
}

\theoremstyle{definition}

\NewTheorem{Definition}{Definitions}
\NewTheorem{Remark}{Remarks}
\NewTheorem{Axiom}{Axioms}
\NewTheorem{Example}{Examples}
\NewTheorem{Observation}{Observations}
\NewTheorem{Convention}{Conventions}
\NewTheorem{Notation}{Notations}
\NewTheorem{Setting}{Settings}
\NewTheorem{Question}{Questions}
\NewTheorem{Answer}{Answers}
\NewTheorem{Conjecture}{Conjectures}
\NewTheorem{Problem}{Problems}
\NewTheorem{Solution}{Solutions}
\NewTheorem{Goal}{Goals}
\NewTheorem{Comment}{Comments}
\NewTheorem{Aim}{Aims}
\NewTheorem{Caution}{Cautions}
\NewTheorem{Exercise}{Exercises}

\theoremstyle{plain}
\NewTheorem{Proposition}{Propositions}
\NewTheorem{Lemma}{Lemmas}
\NewTheorem{Theorem}{Theorems}
\NewTheorem{Corollary}{Corollaries}

\crefname{enumi}{}{}
\Crefname{enumi}{}{}
\creflabelformat{enumi}{(#2#1#3)}
\crefname{enumii}{}{}
\Crefname{enumii}{}{}
\creflabelformat{enumii}{(#2#1#3)}
\crefname{enumiii}{}{}
\Crefname{enumiii}{}{}
\creflabelformat{enumiii}{(#2#1#3)}

\makeatletter
\renewcommand{\p@enumii}{}
\renewcommand{\p@enumiii}{}
\makeatother

\numberwithin{equation}{section}
\crefname{equation}{}{}
\Crefname{equation}{}{}
\creflabelformat{equation}{(#2#1#3)}

\newcommand{\SwapSymbols}[1]{
	\expandafter\let\expandafter\temporarysymbol\csname #1\endcsname
	\expandafter\let\csname #1\expandafter\endcsname\csname var#1\endcsname
	\expandafter\let\csname var#1\endcsname\temporarysymbol
}

\SwapSymbols{epsilon}
\SwapSymbols{phi}
\SwapSymbols{Gamma}
\SwapSymbols{Delta}
\SwapSymbols{Theta}
\SwapSymbols{Lambda}
\SwapSymbols{Xi}
\SwapSymbols{Pi}
\SwapSymbols{Sigma}
\SwapSymbols{Upsilon}
\SwapSymbols{Phi}
\SwapSymbols{Psi}
\SwapSymbols{Omega}

\newcommand{\bbP}{\mathbb{P}}

\newcommand{\bbZ}{\mathbb{Z}}

\newcommand{\cA}{\mathcal{A}}

\newcommand{\cC}{\mathcal{C}}
\newcommand{\cD}{\mathcal{D}}

\newcommand{\cF}{\mathcal{F}}
\newcommand{\cG}{\mathcal{G}}

\newcommand{\cI}{\mathcal{I}}

\newcommand{\cL}{\mathcal{L}}
\newcommand{\cM}{\mathcal{M}}
\newcommand{\cN}{\mathcal{N}}
\newcommand{\cO}{\mathcal{O}}

\newcommand{\cU}{\mathcal{U}}

\newcommand{\cX}{\mathcal{X}}
\newcommand{\cY}{\mathcal{Y}}
\newcommand{\cZ}{\mathcal{Z}}

\newcommand{\km}{\mathfrak{m}}

\let\originalleft\left
\let\originalright\right
\renewcommand{\left}{\mathopen{}\mathclose\bgroup\originalleft}
\renewcommand{\right}{\aftergroup\egroup\originalright}

\newcommand{\into}{\hookrightarrow}

\newcommand{\isoto}{\xrightarrow{\smash{\raisebox{-0.25em}{$\sim$}}}}

\newcommand{\set}[2][]{\mathopen{#1\{}#2\mathclose{#1\}}}

\newcommand{\setwithcondition}[3][]{\mathopen{#1\{}\,#2\mathrel{#1|}#3\,\mathclose{#1\}}}

\newcommand{\op}{\textnormal{op}}

\DeclareMathOperator{\Hom}{Hom}
\DeclareMathOperator{\End}{End}

\DeclareMathOperator{\Ext}{Ext}

\DeclareMathOperator{\Func}{Func}

\DeclareMathOperator{\Mod}{Mod}

\DeclareMathOperator{\MOD}{MOD}

\DeclareMathOperator{\QCoh}{QCoh}

\DeclareMathOperator{\Ann}{Ann}

\DeclareMathOperator{\Spec}{Spec}

\newcommand{\resp}{resp.\ }

\title{Non-exactness of direct products of quasi-coherent sheaves}

\author{Ryo Kanda}
\address[Ryo Kanda]{Department of Mathematics, Graduate School of Science, Osaka University, Toyonaka, Osaka, 560-0043, Japan}
\email{ryo.kanda.math@gmail.com}

\begin{document}

\begin{abstract}
	For a noetherian scheme that has an ample family of invertible sheaves, we prove that direct products in the category of quasi-coherent sheaves are not exact unless the scheme is affine. This result can especially be applied to all quasi-projective schemes over commutative noetherian rings. The main tools of the proof are the Gabriel-Popescu embedding and Roos' characterization of Grothendieck categories satisfying Ab6 and Ab4*.
\end{abstract}

\maketitle
\tableofcontents

\section{Introduction}
\label{21714825}

The class of Grothendieck categories is a large framework that includes
\begin{itemize}
	\item the category $\Mod R$ of right modules over a ring $R$,
	\item the category $\QCoh X$ of quasi-coherent sheaves on a scheme $X$, and
	\item the category of sheaves of abelian groups on a topological space.
\end{itemize}
One of the significant properties of $\Mod R$ for rings $R$ among Grothendieck categories is the exactness of direct products, which is known as Grothendieck's condition Ab4*. This is immediately verified by direct computation, but it is also a consequence of the fact that $\Mod R$ has enough projectives. In general, for a Grothendieck category, the exactness of direct products is equivalent to the category having projective effacements, which is a weak lifting property that resembles the property of having enough projectives. However, it is known that there exists a Grothendieck category that has exact direct products but does not have any nonzero projective objects (see \cref{68007626}). The main source of Grothendieck categories with exact direct products is a pair of a ring and an idempotent ideal of it (\cref{21981963}). Such a pair is used as the \emph{basic setup} of almost ring theory (\cite[2.1.1]{MR2004652}).

For a scheme $X$, it is apparently rare that $\QCoh X$ has exact direct products. Indeed, it is known that direct products in $\QCoh X$ are not exact when $X$ is either
\begin{itemize}
	\item the projective line over a field (\cite{MR2157133}), or
	\item the punctured spectrum of a regular local (commutative noetherian) ring with Krull dimension at least two (\cite{MR0215895});
\end{itemize}
see \cref{17137985,25908527}.

The aim of this paper is to generalize these observations to a wide class of schemes:

\begin{Theorem}[\cref{74600988}]\label{71355712}
	Let $X$ be a divisorial noetherian scheme. Then the following conditions are equivalent:
	\begin{enumerate}
		\item\label{28230044} Direct products in $\QCoh X$ are exact.
		\item\label{41570356} $\QCoh X$ has enough projectives.
		\item\label{28894068} $X$ is an affine scheme.
	\end{enumerate}
\end{Theorem}

A noetherian scheme is called \emph{divisorial} if it admits an ample family of invertible sheaves (\cref{79913009}). Since the exactness of direct products is inherited by closed subschemes, we obtain the following corollary:

\begin{Corollary}[\cref{44529212}]\label{66241253}
	Let $X$ be a scheme that contains a non-affine divisorial noetherian scheme as a closed subscheme. Then direct products in $\QCoh X$ are not exact.
\end{Corollary}

\begin{Remark}\label{20253537}
	Since a divisorial noetherian scheme is a generalization of a noetherian scheme having an ample invertible sheaf, \cref{71355712} can be applied to quasi-projective schemes over commutative noetherian rings. Therefore the aforementioned results of \cite{MR2157133} and \cite{MR0215895} can be derived from \cref{71355712}.
\end{Remark}

\subsection*{Acknowledgments}

The author thanks Shinnosuke Okawa and Akiyoshi Sannai for providing the idea of \cref{44529212}, and the anonymous referee for their valuable comments.

The author was a JSPS Overseas Research Fellow. This work was supported by JSPS KAKENHI Grant Numbers JP17K14164 and JP16H06337.

\section{Gabriel-Popescu embedding and Roos' theorem}
\label{61866504}

\subsection{Preliminaries}
\label{38859488}

\begin{Convention}\label{76880124}\leavevmode
	\begin{enumerate}
		\item\label{79019710} Throughout this paper, we fix a Grothendieck universe. A \emph{small set} is an element of the universe. For objects $X$ and $Y$ in a category $\cC$, the set $\Hom_{\cC}(X,Y)$ is assumed to be small. Colimits and limits always mean those whose index sets are in bijection with small sets. The index set of a generating set is also in bijection with a small set. All rings, schemes, and modules are assumed to be small.
		\item\label{45910523} Since we mainly work on Grothendieck categories (which resemble the category of modules over a ring), coproducts of objects are called direct sums, and products are called direct products. For a family $\set{M_{i}}_{i\in I}$ of objects in a category, its direct sum and direct product are denoted by $\bigoplus_{i\in I}M_{i}$ and $\prod_{i\in I}M_{i}$, respectively. A direct limit (\resp inverse limit) is a colimit (\resp limit) of a direct system (\resp inverse system) indexed by a directed set.
	\end{enumerate}
\end{Convention}

We recall Grothendieck's conditions on exactness of colimits and limits and the definition of a generating set:

\begin{Definition}\label{44044651}\leavevmode
	\begin{enumerate}
		\item\label{90445887} Let $\cA$ be an abelian category that admits direct sums (\resp direct products).
		\begin{enumerate}
			\item\label{47232421} We say that $\cA$ \emph{satisfies Ab4} (\resp \emph{Ab4*}) if direct sums (\resp \emph{direct products}) are exact in $\cA$, that is, for every family of short exact sequences
			\begin{equation*}
				0\to L_{i}\to M_{i}\to N_{i}\to 0\quad\text{($i\in I$)}
			\end{equation*}
			in $\cA$ (where $I$ is in bijection with a small set), the termwise direct sum (\resp direct product)
			\begin{equation*}
				0\to\bigoplus_{i\in I}L_{i}\to\bigoplus_{i\in I}M_{i}\to\bigoplus_{i\in I}N_{i}\to 0
			\end{equation*}
			is again a short exact sequence.
			\item\label{65066674} We say that $\cA$ \emph{satisfies Ab5} (\resp \emph{Ab5*}) if direct limits (\resp \emph{inverse limits}) are exact in $\cA$, that is, for every direct system (\resp inverse system) of short exact sequences in $\cA$, the termwise direct limit (\resp inverse limit) is again exact.
		\end{enumerate}
		\item\label{36313597} Let $\cA$ be an abelian category. A set $\cU$ of objects in $\cA$ that is in bijection with a small set is called a \emph{generating set} if for every nonzero morphism $f\colon X\to Y$ in $\cA$, there exists $U\in\cU$ and a morphism $g\colon U\to X$ such that $fg\neq 0$. An object $U\in\cA$ is called a \emph{generator} if the singleton $\set{U}$ is a generating set.
		\item\label{20908814} An abelian category is called a \emph{Grothendieck category} if it satisfies Ab5 and has a generator.
	\end{enumerate}
\end{Definition}

\begin{Remark}\label{57511646}\leavevmode
	\begin{enumerate}
		\item\label{87000546} If an abelian category admits direct sums (\resp direct products), then it admits colimits (\resp limits) (see \cite[Proposition~2.2.9]{MR2182076}). Direct sums and direct limits (\resp direct products and inverse limits) are always right exact (\resp left exact). So the conditions in \cref{44044651} \cref{90445887} only require left exactness (\resp right exactness).
		\item\label{04629717} For an abelian category with direct sums, Ab5 implies Ab4 (\cite[Corollary~2.8.9]{MR0340375}). See \cite[Theorem~2.8.6]{MR0340375} for conditions equivalent to Ab5.
		\item\label{99834034} It is known that every Grothendieck category admits direct products (\cite[Corollary~3.7.10]{MR0340375}).
		\item\label{57136512} If an abelian category $\cA$ admits direct sums and has a generating set $\set{U_{i}}_{i\in I}$, then the direct sum $\bigoplus_{i\in I}U_{i}$ is a generator in $\cA$ (\cite[Proposition~2.8.2]{MR0340375}).
	\end{enumerate}
\end{Remark}

\begin{Example}\label{41019523}
	Let $R$ be a ring. Then the category $\Mod R$ of right $R$-modules is a Grothendieck category satisfying Ab4*.
\end{Example}

\begin{Example}\label{14163805}
	Let $X$ be a scheme. Then the category $\QCoh X$ of quasi-coherent sheaves on $X$ is an abelian category satisfying Ab5. It was shown by Gabber that $\QCoh X$ has a generator (see \cite[Remarks A.1 and A.2]{MR3780029}). Hence $\QCoh X$ is a Grothendieck category.
\end{Example}

\cref{14163805} implies that $\QCoh X$ for a scheme $X$ admits direct products. However, direct products are not necessarily exact as the following two results show (see also \cite{MR1029695} and \cite[Example~3.10]{MR3922832} for related results):

\begin{Theorem}[{Keller; see \cite[Example~4.9]{MR2157133}}]\label{17137985}
	Let $X:=\bbP^{1}_{k}$ be the projective line over a field $k$. Then $\QCoh X$ does not satisfy Ab4*.
\end{Theorem}

\begin{Theorem}[{Roos \cite[Example~3]{MR0215895}}]\label{25908527}
	Let $R$ be a regular local (commutative noetherian) ring with maximal ideal $\km$. Define $X:=\Spec R\setminus\set{\km}$ as an open subscheme of $\Spec R$. Then $\QCoh X$ satisfies Ab4* if and only if $\dim R\leq 1$.
\end{Theorem}

\begin{Remark}\label{68007626}
	Recall that an abelian category $\cA$ is said to have \emph{enough projectives} if each object in $\cA$ is a quotient object of some projective object. Every Grothendieck category that has enough projectives satisfies Ab4* (the dual of \cite[Corollary~3.2.9]{MR0340375}). The converse does not hold. Indeed, it is shown in \cite[Example~4.2]{MR2197371} that there exists a nonzero Grothendieck category that satisfies Ab4* but has no nonzero projective objects.
	
	It is known that a Grothendieck category satisfies Ab4* if and only if it has \emph{projective effacements} (\cite[Remark~1 in p.~137]{MR0102537}; see also \cite[Corollary~1.4]{MR2197371}), which can be regarded as a weak form of having enough projectives.
\end{Remark}

\begin{Remark}\label{97631379}
	The category of sheaves of abelian groups on a topological space is also a typical example of a Grothendieck category. The exactness of direct products in such a category is characterized in \cite[Corollary~1]{MR0215895} (see also \cite[Theorem~1.7]{MR2197371}).
\end{Remark}

We recall the definitions and basic properties of some classes of subcategories.

\begin{Definition}\label{44486747}
	Let $\cG$ be a Grothendieck category.
	\begin{enumerate}
		\item\label{34448059} A \emph{Serre subcategory} of $\cG$ is a full subcategory of $\cG$ closed under subobjects, quotient objects, and extensions. If $\cX\subset\cG$ is a Serre subcategory, then we have the \emph{quotient category} of $\cG$ by $\cX$, which is denoted by $\cG/\cX$, together with a canonical functor $\cG\to\cG/\cX$ (see \cite[Section~4.3]{MR0340375}).
		\item\label{09184337} A Serre subcategory $\cX\subset\cG$ is called a \emph{localizing subcategory} if the canonical functor $\cG\to\cG/\cX$ admits a right adjoint.
		\item\label{15898117} A localizing subcategory $\cX\subset\cG$ is called a \emph{bilocalizing subcategory} if the canonical functor $\cG\to\cG/\cX$ also admits a left adjoint.
	\end{enumerate}
\end{Definition}

\begin{Remark}\label{67482381}
	If $\cX$ is a localizing subcategory of a Grothendieck category $\cG$, then $\cG/\cX$ is again a Grothendieck category (\cite[Corollary~4.6.2]{MR0340375}). The right adjoint $G\colon\cG/\cX\to\cG$ of the canonical functor $F\colon\cG\to\cG/\cX$ is fully faithful, and thus the counit $FG\to 1_{\cG/\cX}$ is an isomorphism (\cite[Proposition~4.4.3]{MR0340375}).
	
	See \cite[Section~4.3]{MR0340375} or \cite[Theorem~5.11]{MR3351569} for basic properties of the quotient category $\cG/\cX$.
\end{Remark}

\begin{Example}[{See \cite[Example~4.3]{MR3780029}}]\label{34102270}
	Let $X$ be a quasi-separated scheme and let $i\colon U\into X$ be an open immersion from a quasi-compact open subscheme $U$. Then $i_{*}\colon\QCoh U\to\QCoh X$ and its left adjoint $i^{*}\colon\QCoh X\to\QCoh U$ induces an equivalence
	\begin{equation*}
		\frac{\QCoh X}{\cY}\isoto\QCoh U,
	\end{equation*}
	where $\cY\subset\QCoh X$ is the localizing subcategory consisting of all objects $\cM\in\QCoh X$ with $i^{*}\cM=0$.
\end{Example}

\begin{Proposition}\label{29744485}
	Let $\cG$ be a Grothendieck category.
	\begin{enumerate}
		\item\label{20025896} Let $\cX\subset\cG$ be a Serre subcategory. Then the following conditions are equivalent:
		\begin{enumerate}
			\item\label{52200636} $\cX$ is a localizing subcategory.
			\item\label{18345490} $\cX$ is closed under direct sums.
			\item\label{02298799} Every object $M\in\cG$ has a largest subobject belonging to $\cX$.
		\end{enumerate}
		\item\label{01292264} Let $\cX\subset\cG$ be a localizing subcategory. Then the following conditions are equivalent:
		\begin{enumerate}
			\item\label{77486178} $\cX$ is a bilocalizing subcategory.
			\item\label{80123255} $\cX$ is closed under direct products.
			\item\label{56822760} Every object $M\in\cG$ has a largest quotient object belonging to $\cX$, that is, $M$ has a smallest subobject among those $L$ satisfying $M/L\in\cX$.
		\end{enumerate}
	\end{enumerate}
\end{Proposition}

\begin{proof}
	\cref{20025896} \cite[Theorem~4.5.2 and Proposition~4.6.3]{MR0340375}.
	
	\cref{01292264} This can be shown in a similar way to the proof of \cite[Theorem~4.21.1]{MR0340375} for the category of modules over a ring.
\end{proof}

\begin{Definition}\label{60466237}
	Let $\cG$ be a Grothendieck category. A \emph{closed subcategory} of $\cG$ is a full subcategory closed under subobjects, quotient objects, direct sums, and direct products.
\end{Definition}

\begin{Remark}\label{63492107}
	Since the direct sum $\bigoplus_{i\in I}M_{i}$ of objects in a Grothendieck category can be regarded as a subobject of the direct product $\prod_{i\in I}M_{i}$, the condition of being closed under direct sums in \cref{60466237} can be omitted.
	
	By \cref{29744485} \cref{01292264}, a full subcategory of a Grothendieck category is bilocalizing if and only if it is localizing and closed.
\end{Remark}

\begin{Proposition}\label{17471113}
	Let $\cG$ be a Grothendieck category and let $\cC\subset\cG$ be a full subcategory closed under subobjects and quotient objects. Then the following conditions are equivalent:
	\begin{enumerate}
		\item\label{74970637} $\cC$ is a closed subcategory.
		\item\label{75321870} Every object $M\in\cG$ has a largest quotient object belonging to $\cC$.
	\end{enumerate}
\end{Proposition}

\begin{proof}
	\cite[Proposition~11.2]{MR3452186}.
\end{proof}

\begin{Remark}\label{57307712}
	For a ring $R$, there exists a bijective correspondence between the two-sided ideals of $R$ and the closed subcategories of $\Mod R$ that sends each $I$ to $\Mod(R/I)$. The bilocalizing subcategories correspond to the idempotent ideals. See \cite[Theorem~11.3 and Proposition~12.6]{MR3452186} for more details.
	
	For a scheme $X$, it is known that the closed subcategories of $\QCoh X$ are in bijection with the closed subschemes of $X$ provided that one of the following conditions holds:
	\begin{itemize}
		\item $X$ is locally noetherian (\cite[Theorem~1.3]{MR3452186}).
		\item $X$ is separated (\cite[Proposition~A.5]{MR3780029}).
	\end{itemize}
\end{Remark}

\subsection{Gabriel-Popescu embedding}
\label{38009298}

A generalization of the Gabriel-Popescu embedding is one of the main tools to prove \cref{71355712}. First we recall the original version:

\begin{Theorem}[Gabriel-Popescu embedding \cite{MR0166241}]\label{85724040}
	Let $\cG$ be a Grothendieck category and let $U\in\cG$ be a generator. Then the functor $\Hom_{\cG}(U,-)\colon\cG\to\Mod\End_{\cG}(U)$ induces an equivalence
	\begin{equation*}
		\cG\isoto\frac{\Mod\End_{\cG}(U)}{\cX},
	\end{equation*}
	where $\cX\subset\Mod\End_{\cG}(U)$ is the localizing subcategory consisting of all $M\in\Mod\End_{\cG}(U)$ annihilated by the left adjoint of $\Hom_{\cG}(U,-)$.
\end{Theorem}

If we have a generating set $\set{U_{i}}_{i\in I}$ in a Grothendieck category $\cG$, then the direct sum $\bigoplus_{i\in I}U_{i}$ is a generator in $\cG$ and we can apply the Gabriel-Popescu embedding. On the other hand, there is a generalized version of the embedding that respects the structure of the given generating set (\cref{70192993}). To state the result, we recall some basic facts on rings that do not necessarily have an identity element.

\begin{Definition}\label{05527173}
	Let $R$ be a ring not necessarily with identity.
	\begin{enumerate}
		\item\label{45126180} A \emph{complete set of orthogonal idempotents} in $R$ is a set of idempotents $\set{e_{i}}_{i\in I}\subset R$ such that
		\begin{itemize}
			\item $\set{e_{i}}_{i\in I}$ is orthogonal, that is, $e_{i}e_{j}=0$ for $i\neq j$, and
			\item $R=\bigoplus_{i,j\in I}e_{i}Re_{j}$ (or equivalently, $R=\bigoplus_{i\in I}e_{i}R=\bigoplus_{i\in I}Re_{i}$).
		\end{itemize}
		We say that \emph{$R$ has enough idempotents} if it admits a complete set of orthogonal idempotents.
		\item\label{47808356} Suppose that $R$ has enough idempotents. The category of all right $R$-modules is denoted by $\Mod R$. Define $\MOD R\subset\Mod R$ to be the full subcategory consisting of right $R$-modules $M$ with $MR=M$.
	\end{enumerate}
\end{Definition}

\begin{Remark}\label{91501107}
	Let $R$ be a ring not necessarily with identity and let $\set{e_{i}}_{i\in I}$ be a complete set of orthogonal idempotents.
	\begin{enumerate}
		\item\label{84098455} For every $M\in\Mod R$, the condition $M=MR$ is equivalent to $M=\bigoplus_{i\in I}Me_{i}$.
		\item\label{86849709} The \emph{Dorroh overring} $R^{*}$ of $R$ is $\bbZ\times R$ as an abelian group and has multiplication
		\begin{equation*}
			(n_{1},r_{1})\cdot (n_{2},r_{2}):=(n_{1}n_{2},n_{1}r_{2}+n_{2}r_{1}+r_{1}r_{2}).
		\end{equation*}
		The Dorroh overring $R^{*}$ is a ring with identity $(1,0)$, and $R$ is identified with $0\times R\subset R^{*}$ (see \cite[Sections 1.5 and 6.3]{MR1144522}). Since the forgetful functor $\Mod R^{*}\to \Mod R$ is an equivalence, the category $\Mod R$ is a Grothendieck category satisfying Ab4*, and colimits and limits in $\Mod R$ can be computed in the same way as the category of right modules over a ring with identity.
		
		\item\label{81049440} $\MOD R\subset\Mod R$ is closed under subobjects, quotient objects, extensions, and direct sums (\cite[Proposition~0.1]{MR2052174}) and $R$ is a projective generator in $\MOD R$. Hence $\MOD R$ is also a Grothendieck category satisfying Ab4*, but limits in $\MOD R$ are different from those computed in $\Mod R$ in general.
	\end{enumerate}
\end{Remark}

\begin{Theorem}[{N{\u{a}}st{\u{a}}sescu and Chite{\c{s}} \cite[Theorem~2.1]{MR2785804}}]\label{70192993}
	Let $\cG$ be a Grothendieck category and let $\set{U_{i}}_{i\in I}$ be a generating set in $\cG$. Then $R:=\bigoplus_{i,j\in I}\Hom_{\cG}(U_{i},U_{j})$ is a ring with enough idempotents, and the functor $G:=\bigoplus_{i\in I}\Hom_{\cG}(U_{i},-)\colon\cG\to\MOD R$ induces an equivalence
	\begin{equation*}
		\cG\isoto\frac{\MOD R}{\cX},
	\end{equation*}
	where $\cX\subset\MOD R$ is the localizing subcategory consisting of all $M\in\MOD R$ annihilated by the left adjoint of $G$.
\end{Theorem}

\begin{Remark}\label{11368333}
	For $\bbZ$-linear categories $\cC$ and $\cD$, let $\Func_{\bbZ}(\cC,\cD)$ denote the category of $\bbZ$-functors from $\cC$ to $\cD$. In the setting of \cref{70192993}, we have an equivalence $\Func_{\bbZ}(\cU^{\op},\Mod\bbZ)\isoto\MOD R$ given by $F\mapsto\bigoplus_{i\in I}F(U_{i})$, where $\cU:=\set{U_{i}}_{i\in I}$ is regarded as a full subcategory of $\cG$. Hence \cref{70192993} can be interpreted in terms of the functor category.
\end{Remark}

\begin{Theorem}[{Prest \cite[Theorem~1.1]{MR585223}; see also \cite[Corollary~2.5]{MR2785804}}]\label{28835132}
	Let $\cG$ be a Grothendieck category and let $\cU$ be a generating set in $\cG$. Then the functor $G\colon\cG\to\Func_{\bbZ}(\cU^{\op},\Mod\bbZ)$ defined by $M\mapsto\Hom_{\cG}(-,M)$ induces an equivalence
	\begin{equation*}
		\cG\isoto\frac{\Func_{\bbZ}(\cU^{\op},\Mod\bbZ)}{\cX},
	\end{equation*}
	where $\cX\subset\Func_{\bbZ}(\cU^{\op},\Mod\bbZ)$ is the localizing subcategory consisting of all objects annihilated by the left adjoint of $G$.
\end{Theorem}

\subsection{Roos' theorem}
\label{80551652}

A generalization of Roos' theorem (\cref{08099896}) is another main ingredient of the proof of \cref{71355712}. To state the result, we recall Grothendieck's condition Ab6:

\begin{Definition}\label{00868717}
	We say that an abelian category $\cA$ with direct sums \emph{satisfies Ab6} if the following assertion holds for every object $M\in\cA$: For every family $\set{\set{L_{i}^{j}}_{i\in I_{j}}}_{j\in J}$ of directed sets of subobjects of $M$ with respect to inclusion, we have
	\begin{equation*}
		\bigcap_{j\in J}\Bigg(\sum_{i\in I_{j}}L_{i}^{j}\Bigg)=\sum_{(i_{j})_{j\in J}\in\prod_{j\in J}I_{j}}\Bigg(\bigcap_{j\in J}L_{i_{j}}^{j}\Bigg).
	\end{equation*}
\end{Definition}

\begin{Remark}\label{41186065}
	For an abelian category with direct sums, Ab6 implies Ab5 (\cite[Corollary~2.8.13]{MR0340375}).
\end{Remark}

\begin{Remark}\label{02181172}
	Roos \cite[Theorem~1]{MR0217145} showed that condition Ab6 has the following characterization (see also \cite[Exercise~3.5.7]{MR0340375}): Let $\cG$ be a Grothendieck category. A subobject $L$ of an object $M\in\cG$ is said to be \emph{of finite type relative to $M$} if for every directed set $\set{L_{i}}_{i\in I}$ of subobjects of $M$ with respect to inclusion satisfying $\sum_{i\in I}L_{i}=M$, the directed set $\set{L_{i}\cap L}_{i\in I}$ eventually stabilizes, that is, $L_{i}\cap L=L$ for some $i\in I$. Then $\cG$ satisfies Ab6 if and only if every object $M\in\cG$ is the sum of all subobjects of finite type relative to $M$.
\end{Remark}

\begin{Remark}\label{64547726}
	Grothendieck categories that we encounter in practice often satisfy Ab6:
	\begin{enumerate}
		\item\label{25819375} Let $\cU$ be a small $\bbZ$-category and let $\cG$ be a Grothendieck category satisfying Ab6. Then $\Func_{\bbZ}(\cU,\cG)$ is a Grothendieck category satisfying Ab6 (\cite[Theorem~3.4.2]{MR0340375}). In particular, for every ring $R$, $\Mod R$ is a Grothendieck category satisfying Ab6.
		\item\label{93398935} A Grothendieck category is called \emph{locally noetherian} if admits a generating set consisting of noetherian objects. Every locally noetherian Grothendieck category satisfies Ab6. This follows from \cref{02181172} since all noetherian subobjects of an object $M$ are of finite type relative to $M$.
		\item\label{38921065} Let $\cG$ be a Grothendieck category satisfying Ab6 (\resp Ab4*) and let $\cX\subset\cG$ be a bilocalizing subcategory. Then $\cG/\cX$ is a Grothendieck category satisfying Ab6 (\resp Ab4*). This follows because the canonical functor $\cG\to\cG/\cX$ preserves all colimits and limits.
	\end{enumerate}
\end{Remark}

\begin{Remark}\label{21981963}
	Let $R$ be a ring and let $I\subset R$ be an idempotent ideal. Then, as in \cref{57307712}, $\Mod(R/I)$ is a bilocalizing subcategory of $\Mod R$. Thus the quotient category of $\Mod R$ by $\Mod(R/I)$ is a Grothendieck category satisfying Ab6 and Ab4* by \cref{64547726}.
\end{Remark}

Roos' theorem shows that all Grothendieck categories satisfying Ab6 and Ab4* arise in the way of \cref{21981963}:

\begin{Theorem}[{Roos \cite[Theorem~1]{MR0190207}}]\label{11569500}
	Let $\cG$ be a Grothendieck category and let $U\in\cG$ be a generator. Define the localizing subcategory $\cX\subset\Mod\End_{\cG}(U)$ as in \cref{85724040}. Then the following conditions are equivalent:
	\begin{enumerate}
		\item $\cG$ satisfies Ab6 and Ab4*.
		\item $\cX$ is closed under direct products, that is, $\cX\subset\Mod\End_{\cG}(U)$ is a bilocalizing subcategory.
	\end{enumerate}
\end{Theorem}

Roos' theorem can be generalized so that it fits into the setting of the generalized Gabriel-Popescu embedding:

\begin{Theorem}\label{08099896}
	Let $\cG$ be a Grothendieck category and let $\set{U_{i}}_{i\in I}$ be a generating set in $\cG$. Let $R:=\bigoplus_{i,j\in I}\Hom_{\cG}(U_{i},U_{j})$. Define the localizing subcategory $\cX\subset\MOD R$ as in \cref{70192993}. Then the following conditions are equivalent:
	\begin{enumerate}
		\item $\cG$ satisfies Ab6 and Ab4*.
		\item $\cX$ is closed under direct products, that is, $\cX\subset\MOD R$ is a bilocalizing subcategory.
	\end{enumerate}
\end{Theorem}

\begin{proof}
	The proof of \cref{11569500} written in \cite[Theorem~4.21.6]{MR0340375} also works in this setting. The proof is modified as follows:
	\begin{enumerate}[label=(\alph*),ref=\alph*]
		\item\label{28558139} Use the generating set $\set{U_{i}}_{i\in I}$ instead of the generator $U$. Use $\MOD R$ instead of $\Mod A$.
		\item\label{20753853} Define $S$ and $\cF$ to be $G$ and $\cX$ in \cref{70192993}, respectively.
		\item\label{49454448} In the conclusion of \cite[Lemma~4.21.3]{MR0340375}, $f$ should run over all elements of $M$ that are homogeneous in the sense that each of them belongs to $\Hom_{\cG}(U_{i},X)$ for some $i\in I$.
		\item\label{66599313} In \cite[Lemma~4.21.4]{MR0340375}, define $X'$ to be $\bigoplus_{i\in I}U_{i}^{\oplus\Hom_{\cG}(U_{i},X)}$.
		\item In \cite[Lemma~4.21.5]{MR0340375}, define $R(X)$ to be the submodule consisting of all finite sums of elements $f\colon U_{i}\to X$ for various $i$ satisfying the same property with $U$ replaced by $U_{i}$. The conclusion of \cite[Lemma~4.21.5]{MR0340375} is modified in the same way as \cref{49454448}.\qedhere
	\end{enumerate}
\end{proof}

\cref{08099896} can also be stated in terms of a functor category:

\begin{Corollary}\label{19244512}
	Let $\cG$ be a Grothendieck category and let $\cU$ be a generating set in $\cG$. Define the localizing subcategory $\cX\subset\Func_{\bbZ}(\cU^{\op},\Mod\bbZ)$ as in \cref{28835132}. Then the following conditions are equivalent:
	\begin{enumerate}
		\item $\cG$ satisfies Ab6 and Ab4*.
		\item $\cX$ is closed under direct products, that is, $\cX\subset\Func_{\bbZ}(\cU^{\op},\Mod\bbZ)$ is a bilocalizing subcategory.
	\end{enumerate}
\end{Corollary}

\begin{proof}
	This is immediate from \cref{08099896} in view of \cref{11368333}.
\end{proof}

\section{Divisorial noetherian schemes}
\label{23512597}

In this section, we prove the main results. Whenever we consider a scheme $X$, unadorned tensor products are tensor products of quasi-coherent sheaves on $X$. The structure sheaf of $X$ is denoted by $\cO_{X}$, and $\Gamma(X,-)$ is the global section functor.

We recall the definition of a divisorial scheme:

\begin{Definition}[{\cite[Proposition~1.1]{MR1970862}; see also \cite[Definition~3.1]{MR0219545}}]\label{79913009}
	Let $X$ be a quasi-compact and quasi-separated scheme.
	\begin{enumerate}
		\item\label{74348417} A finite family $\set{\cL_{1},\ldots,\cL_{r}}$ of invertible sheaves on $X$ is called an \emph{ample family} if the set
		\begin{equation*}
			\setwithcondition{X_{s}}{\textnormal{$s\in\Gamma(X,\cL_{1}^{\otimes d_{1}}\otimes\cdots\otimes\cL_{r}^{\otimes d_{r}})$, $d_{1},\ldots,d_{r}\geq 0$ are integers}}
		\end{equation*}
		is an open basis of $X$, where $X_{s}\subset X$ is the open subset consisting of all $x\in X$ such that $s_{x}$ does not belong the unique maximal ideal of $(\cL_{1}^{\otimes d_{1}}\otimes\cdots\otimes\cL_{r}^{\otimes d_{r}})_{x}$.
		\item\label{67256126} $X$ is called \emph{divisorial} if it admits an ample family of invertible sheaves.
	\end{enumerate}
\end{Definition}

\begin{Remark}\label{15943873}
	An ample family of invertible sheaves is a generalization of an ample invertible sheaf (see \cite[Section~4.5]{MR0217084}). In particular, every quasi-projective scheme over a commutative noetherian ring is divisorial.
\end{Remark}

The following fact is essential for our proof:

\begin{Proposition}[{\cite[Theorem~3.3]{MR0219545}}]\label{23753093}
	Let $X$ be a divisorial noetherian scheme. Then every coherent sheaf on $X$ is isomorphic to a quotient of a direct sum of invertible sheaves.
\end{Proposition}

\begin{Remark}\label{96863411}
	For a noetherian scheme $X$, the category $\QCoh X$ is a locally noetherian Grothendieck category (\cite[Theorem~1 in p.~443]{MR0232821}). An object in $\QCoh X$ is noetherian if and only if it is a coherent sheaf on $X$.
	
	Hence \cref{23753093} implies that if $X$ is a divisorial noetherian scheme, then $\QCoh X$ has a generating set consisting of invertible sheaves.
\end{Remark}

\begin{Setting}\label{48175114}
	In the rest of this section, let $X$ be a divisorial noetherian scheme. We use the following notations:
	\begin{enumerate}
		\item\label{15319342} Fix a generating set $\set{\cL_{\lambda}}_{\lambda\in\Lambda}$ in $\QCoh X$ consisting of invertible sheaves (see \cref{96863411}).
		\item\label{17530454} Let $I$ be the free abelian group generated by $\Lambda$, that is, $I=\bigoplus_{\lambda\in\Lambda}\bbZ\lambda$. For each $i=\sum_{j}n_{j}\lambda_{j}\in I$, where $n_{j}\in\bbZ$ and $\lambda_{j}\in\Lambda$, define the invertible sheaf
		\begin{equation*}
			\cL_{i}:=\bigotimes_{j}\cL_{\lambda_{j}}^{\otimes n_{j}}.
		\end{equation*}
		Then $\set{\cL_{i}}_{i\in I}=\set{\cL_{-i}}_{i\in I}$ is also a generating set in $\QCoh X$. Define the ring not necessarily with identity
		\begin{equation*}
			R:=\bigoplus_{i,j\in I}\Hom_{X}(\cL_{-i},\cL_{-j}).
		\end{equation*}
		Let $e_{i}\in\Hom_{X}(\cL_{-i},\cL_{-i})$ be the identity morphism for every $i\in I$. Then $\set{e_{i}}_{i\in I}$ is a complete set of orthogonal idempotents of $R$.
		\item\label{97146589} Define the $I$-graded ring
		\begin{equation*}
			S:=\bigoplus_{i\in I}\Gamma(X,\cL_{i})=\bigoplus_{i\in I}\Hom_{X}(\cO_{X},\cL_{i})
		\end{equation*}
		with the following multiplication: For each $f\in\Hom_{X}(\cO_{X},\cL_{i})$ and $g\in\Hom_{X}(\cO_{X},\cL_{j})$, $gf\in\Hom_{X}(\cO_{X},\cL_{i+j})$ is defined to be the composite
		\begin{equation*}
			\begin{tikzcd}
				\cO_{X}\ar[r,"f"] & \cL_{i}\ar[d,"\wr"] & \cL_{i+j} \\
				& \cL_{i}\otimes\cO_{X}\ar[r,"\cL_{i}\otimes g"'] & \cL_{i}\otimes\cL_{j}\rlap{,}\ar[u,"\wr"]
			\end{tikzcd}
		\end{equation*}
		where the isomorphisms are the canonical ones. It is straightforward to see that $S$ is a commutative ring (\emph{with} identity). Denote by $\Mod^{I}S$ the category of $I$-graded $S$-modules whose morphisms are homogeneous $S$-homomorphisms of degree $0$. For an object $M\in\Mod^{I}S$ and $j\in I$, define the degree shift $M(j)\in\Mod^{I}S$ to be the same $S$-module with new grading $M(j)_{i}=M_{i+j}$. This defines the equivalence
		\begin{equation*}
			(j)\colon\Mod^{I}S\isoto\Mod^{I}S.
		\end{equation*}
	\end{enumerate}
\end{Setting}

\begin{Remark}\label{37549406}
	The \emph{$I$-algebra} associated to the $I$-graded ring $S$ is the ring $A$ not necessarily with identity defined by
	\begin{equation*}
			A:=\bigoplus_{i,j\in I}A_{i,j},\quad\text{where}\quad A_{i,j}:=S_{j-i},
		\end{equation*}
	The multiplication $A_{i,j}\times A_{j',k}\to A$ is given by that of $S$ for $j=j'$ and the zero map for $j\neq j'$ (see the last paragraph in \cite[p.~3988]{MR2836401}). There is an isomorphism $A\isoto R$ of rings not necessarily with identities given by
	\begin{equation*}
		-\otimes\cL_{-j}\colon A_{i,j}=\Hom_{X}(\cO_{X},\cL_{j-i})\isoto\Hom_{X}(\cL_{-j},\cL_{-i})=e_{i}Re_{j}.
	\end{equation*}
	
	There is an equivalence $\MOD R\isoto\Mod^{I}S$ that sends each $M\in\MOD R$ to $\bigoplus_{i\in I}Me_{i}$, where the $S$-action $Me_{i}\times S_{j}\to Me_{i+j}$ is induced from the $A$-action $Me_{i}\times A_{i,i+j}\to Me_{i+j}$, or the $R$-action $Me_{i}\times\Hom_{X}(\cL_{-i-j},\cL_{-i})\to Me_{i+j}$.
\end{Remark}

\begin{Lemma}\label{92619170}
	Assume that $\QCoh X$ satisfies Ab4*. Then there exists an equivalence
	\begin{equation*}
		\QCoh X\isoto\frac{\Mod^{I}S}{\cY},
	\end{equation*}
	where $\cY\subset\Mod^{I}S$ is a bilocalizing subcategory closed under degree shifts, that sends $\cO_{X}\in\QCoh X$ to an object isomorphic to the image of $S\in\Mod^{I}S$ by the canonical functor to $(\Mod^{I}S)/\cY$.
\end{Lemma}

\begin{proof}
	$\QCoh X$ satisfies Ab6 since it is locally noetherian (\cref{64547726} \cref{93398935}). Applying \cref{08099896} to the generating set $\set{\cL_{-i}}_{i\in I}$, we deduce that the functor $\bigoplus_{i\in I}\Hom_{X}(\cL_{-i},-)\colon\QCoh X\to\MOD R$ induces an equivalence
	\begin{equation*}
		\QCoh X\isoto\frac{\MOD R}{\cY'}
	\end{equation*}
	for some bilocalizing subcategory $\cY'\subset\MOD R$. The equivalence $\MOD R\isoto\Mod^{I}S$ in \cref{37549406} induces an equivalence
	\begin{equation*}
		\frac{\MOD R}{\cY'}\isoto\frac{\Mod^{I}S}{\cY}
	\end{equation*}
	for some bilocalizing subcategory $\cY\subset\Mod^{I}S$.
	
	Denote by $G$ the composite
	\begin{equation*}
		\QCoh X\to\MOD R\isoto\Mod^{I}S
	\end{equation*}
	and let $F$ be the left adjoint of $G$. Then $G(\cO_{X})=\bigoplus_{i\in I}\Hom_{X}(\cL_{-i},\cO_{X})$, which is isomorphic to $S$ via
	\begin{equation*}
		-\otimes\cL_{i}\colon\Hom_{X}(\cL_{-i},\cO_{X})\isoto\Hom_{X}(\cO_{X},\cL_{i}).
	\end{equation*}
	
	Let $j\in I$. For every object $\cM\in\QCoh X$,
	\begin{equation*}
		G(\cM\otimes\cL_{j})_{i}=\Hom_{X}(\cL_{-i},\cM\otimes\cL_{j})\cong\Hom_{X}(\cL_{-i-j},\cM)=G(\cM)_{i+j}=G(\cM)(j)_{i},
	\end{equation*}
	and it is straightforward to see that this gives an isomorphism $G(\cM\otimes\cL_{j})\cong G(\cM)(j)$ that is functorial in $\cM$. Thus the diagram
	\begin{equation*}
		\begin{tikzcd}
			\QCoh X\ar[d,"-\otimes\cL_{j}"',"\wr"]\ar[r,"G"] & \Mod^{I}S\ar[d,"(j)","\wr"'] \\
			\QCoh X\ar[r,"G"'] & \Mod^{I}S
		\end{tikzcd}
	\end{equation*}
	commutes up to isomorphism. The adjoint property implies that the diagram
	\begin{equation*}
		\begin{tikzcd}
			\QCoh X\ar[d,"-\otimes\cL_{j}"',"\wr"] & \Mod^{I}S\ar[d,"(j)","\wr"']\ar[l,"F"'] \\
			\QCoh X & \Mod^{I}S\ar[l,"F"]
		\end{tikzcd}
	\end{equation*}
	also commutes up to isomorphism. Since $\cY$ consists of all objects in $\Mod^{I}S$ annihilated by $F$, it is closed under degree shifts.
\end{proof}

\begin{Setting}\label{62687661}
	In the subsequent lemmas, we assume that $\QCoh X$ satisfies Ab4*, and use the following notations in addition to \cref{48175114}:
	\begin{enumerate}
		\item\label{40173092} Let $\cY\subset\Mod^{I}S$ be the bilocalizing subcategory closed under degree shifts obtained in \cref{92619170}. Let $F\colon\Mod^{I}S\to\QCoh X$ be the composite
		\begin{equation*}
			\Mod^{I}S\to\frac{\Mod^{I}S}{\cY}\isoto\QCoh X
		\end{equation*}
		of the canonical functor and the equivalence obtained in \cref{92619170}. Let $G$ be its right adjoint. These are the same functors as those appeared in the proof of \cref{92619170}. Note that $G(\cO_{X})\cong S$ by the proof of \cref{92619170}, and $F(S)\cong FG(\cO_{X})\cong\cO_{X}$ is a noetherian object by \cref{67482381,96863411}.
		\item\label{76549029} Define $\cZ\subset\Mod^{I}S$ to be the full subcategory consisting of all objects $M\in\Mod^{I}S$ such that none of the nonzero subquotients of $M$ belong to $\cY$.
	\end{enumerate}
\end{Setting}

We will show that $\cY=0$ and $\cZ=\Mod^{I}S$ in \cref{50382042}. So \cref{26543493,09447381,27915695,01483581} are only used to prove \cref{50382042} and they will eventually become trivial.

\begin{Lemma}\label{26543493}
	Let $M\in\Mod^{I}S$ be an object that belongs to $\cY$. Then every $N\in\Mod^{I}S$ satisfying $N\Ann_{S}(M)=0$ belongs to $\cY$.
\end{Lemma}

\begin{proof}
	Since $\Ann_{S}(M)=\bigcap_{x}\Ann_{S}(x)$, where $x$ runs over all homogeneous elements of $M$, we have the canonical monomorphism
	\begin{equation*}
		\frac{S}{\Ann_{S}(M)}\to\prod_{x}\frac{S}{\Ann_{S}(x)},
	\end{equation*}
	and $S/\Ann_{S}(x)\cong (xS)(\deg x)\subset M(\deg x)\in\cY$. Hence $S/\Ann_{S}(M)$ belongs to $\cY$. The condition $N\Ann_{S}(M)=0$ implies that $N$ is a quotient of a direct sum of copies of $S/\Ann_{S}(M)$. Therefore $N\in\cY$.
\end{proof}

\begin{Lemma}\label{09447381}\leavevmode
	\begin{enumerate}
		\item\label{41600721} $\cZ\subset\Mod^{I}S$ is a localizing subcategory closed under degree shifts.
		\item\label{01957179} Let $M\in\Mod^{I}S$ be a noetherian object that belongs to $\cZ$. Then every $N\in\Mod^{I}S$ satisfying $N\Ann_{S}(M)=0$ belongs to $\cZ$.
	\end{enumerate}
\end{Lemma}

\begin{proof}
	\cref{41600721} It is obvious that $\cZ$ is closed under subobjects and quotient objects. By \cite[Proposition~2.4 (4)]{MR2964615}, every nonzero subquotient of an extension of two objects in $\cZ$ is a nonzero extension of subquotients of objects in $\cZ$, which does not belong to $\cY$. Hence $\cZ$ is closed under extensions. In a similar way to the proof of \cite[Proposition~2.12 (1)]{MR3351569}, we deduce that $\cZ$ is also closed under direct sums. Since $\cY$ is closed under degree shifts, $\cZ$ is also closed under degree shifts.
	
	\cref{01957179} Since $M$ is noetherian, there are a finite number of elements $x_{1},\ldots,x_{n}\in S$ such that $\Ann_{S}(M)=\bigcap_{i=1}^{n}\Ann_{S}(x_{i})$. Therefore the claim can be shown similarly to \cref{26543493}.
\end{proof}

\begin{Remark}\label{52819905}
	If an object $M\in\Mod^{I}S$ has an ascending chain of subobjects $L_{0}\subset L_{1}\subset\cdots$ such that $L_{i+1}/L_{i}\notin\cY$ for all $i\geq 0$, then we obtain the strictly ascending chain $F(L_{0})\subsetneq F(L_{1})\subsetneq\cdots$ of subobjects of $F(M)$ since $F(L_{i+1})/F(L_{i})\cong F(L_{i+1}/L_{i})\neq 0$. Hence, if $F(M)\in\QCoh X$ is noetherian, then there are no such chains of subobjects of $M$.
	
	In particular, if $M\in\cZ$, then $M$ is noetherian if and only if $F(M)$ is noetherian (see \cite[Lemma~5.8.3]{MR0340375} for the ``only if'' part).
\end{Remark}

\begin{Lemma}\label{27915695}
	Let
	\begin{equation*}
		0\to L\to M\to N\to 0
	\end{equation*}
	be a short exact sequence in $\Mod^{I}S$ such that $F(M)\in\QCoh X$ is noetherian and one of the following conditions is satisfied:
	\begin{enumerate}
		\item\label{40504291} $L\in\cY$ and $N\in\cZ$.
		\item\label{83462481} $L\in\cZ$ and $N\in\cY$.
	\end{enumerate}
	Then the exact sequence splits.
\end{Lemma}

\begin{proof}
	Since $S$ is a commutative ring, we have
	\begin{equation*}
		M\Ann_{S}(L)\Ann_{S}(N)=M\Ann_{S}(N)\Ann_{S}(L)=0.
	\end{equation*}
	Assume \cref{40504291}. Then $N$ is noetherian by \cref{52819905}, and \cref{09447381} implies $L':=M\Ann_{S}(L)\in\cZ$. Since $N':=M/L'$ is annihilated by $\Ann_{S}(L)$, \cref{26543493} implies $N'\in\cY$.
	
	Let $K$ be the kernel of the composite $L\to M\to N'$. Then the composite $K\into L\to M$ factors through some morphism $K\to L'$. Since $L\in\cY$ and $\cY$ is closed under subobjects, $K\in\cY$. Thus the only morphism from $K$ to $L'\in\cZ$ is zero. This means $K=0$. The dual argument shows that the cokernel of the composite $L\to M\to N'$ is also zero. Therefore it is an isomorphism, and the given exact sequence splits.
	
	The proof for \cref{83462481} is similar.
\end{proof}

\begin{Lemma}\label{01483581}
	Let $\cC\subset\Mod^{I}S$ be the collection of objects $H$ such that
	\begin{enumerate}
		\item\label{30191978} no nonzero subobjects of $H$ belong to $\cY$, and
		\item\label{71369553} no nonzero subobjects of $H$ belong to $\cZ$.
	\end{enumerate}
	If $M\in\Mod^{I}S$ does not have any nonzero subquotient that belongs to $\cC$ and $F(M)\in\QCoh X$ is a noetherian object, then $M$ is a direct sum of an object in $\cY$ and an object in $\cZ$.
\end{Lemma}

\begin{proof}
	Assume that $M$ satisfies the assumption but it is not a direct sum of an object in $\cY$ and an object in $\cZ$. Since $F(M)$ is noetherian, we can assume that for every nonzero subobject $L\subset M$ with $F(L)\neq 0$, the quotient $M/L$ is a direct sum of an object in $\cY$ and an object in $\cZ$. Indeed, if it is not the case, we can replace $M$ by $M/L$ since $M/L$ satisfies the same assumption. This procedure eventually terminates due to \cref{52819905}.
	
	Take the largest subobject $L\subset M$ belonging to $\cY$ using \cref{29744485} \cref{20025896}. Assume $L=0$. Then $M$ satisfies \cref{30191978}, and hence it does not satisfy \cref{71369553}. $M$ has a nonzero subobject $N\subset M$ that belongs to $\cZ$. Since $\cZ\subset\Mod^{I}S$ is a localizing subcategory by \cref{09447381} \cref{41600721}, $N$ can be taken to be largest among the subobjects belonging to $\cZ$. Since $N\neq 0$, we have $M/N\cong M_{1}\oplus M_{2}$ for some $M_{1}\in\cY$ and $M_{2}\in\cZ$, and the maximality of $N$ implies that $M_{2}=0$. By \cref{27915695}, the short exact sequence
	\begin{equation*}
		0\to N\to M\to M_{1}\to 0
	\end{equation*}
	splits. This contradicts the assumption that $L=0$.
	
	Assume $L\neq 0$. Then the argument for $L=0$ shows that $M/L\cong M_{1}\oplus M_{2}$ for some $M_{1}\in\cY$ and $M_{2}\in\cZ$. By the maximality of $L$, we have $M_{1}=0$. By \cref{27915695}, the short exact sequence
	\begin{equation*}
		0\to L\to M\to M_{2}\to 0
	\end{equation*}
	splits. This is again a contradiction.
\end{proof}

\begin{Lemma}\label{50382042}
	$\cY=0$ and $\cZ=\Mod^{I}S$.
\end{Lemma}

\begin{proof}
	Assume that $S\in\Mod^{I}S$ has a nonzero subquotient $H$ that belongs to the collection $\cC$ defined in \cref{01483581}. If there is a nonzero subobject $L\subset H$ such that $H/L$ has a nonzero subquotient belonging to $\cC$, then replace $H$ by $H/L$. Since $F(S)\cong\cO_{X}$ is noetherian and $L$ does not belong to $\cY$, this procedure eventually terminates by \cref{52819905}. Thus we can assume that for every nonzero subobject $L\subset H$, the quotient $H/L$ does not have any nonzero subquotient belonging to $\cC$.
	
	Since $\cY\subset\Mod^{I}S$ is a bilocalizing subcategory, $H$ has the smallest subobject $H'\subset H$ among those satisfying $H/H'\in\cY$ by \cref{29744485} \cref{01292264}. Since $H'\neq 0$ and $H'$ also belongs to $\cC$, we can assume that no nonzero quotient object of $H$ belong to $\cY$ by replacing $H$ by $H'$.
	
	By property \cref{71369553} in the definition of $\cC$, there exist subobjects $L\subsetneq L'\subset H$ such that $L'/L\in\cY$. Then $H/L$ meets the requirement on $M$ in \cref{01483581}. Hence $H/L\cong M_{1}\oplus M_{2}$ for some $M_{1}\in\cY$ and $M_{2}\in\cZ$. Since $H/L$ has a nonzero subobject $L'/L\in\cY$, the direct summand $M_{1}$ is nonzero. This contradicts to that $H$ has no nonzero quotient object that belongs to $\cY$.
	
	Therefore $S\in\Mod^{I}S$ does not have any nonzero subquotient that belongs to $\cC$. Again by \cref{01483581}, $S\cong N_{1}\oplus N_{2}$ for some $N_{1}\in\cY$ and $N_{2}\in\cZ$. Since $F(N_{1})=0$, we have
	\begin{equation*}
		\Hom_{S}(N_{1},S)\cong\Hom_{S}(N_{1},G(\cO_{X}))\cong\Hom_{S}(F(N_{1}),\cO_{X})=0.
	\end{equation*}
	Hence $N_{1}=0$, and $S=N_{2}\in\cZ$. Since $\set{S(i)}_{i\in I}$ is a generating set in $\Mod^{I}S$ and $\cZ\subset\Mod^{I}S$ is a localizing subcategory closed under degree shifts by \cref{09447381} \cref{41600721}, we obtain $\cZ=\Mod^{I}S$. This implies that $\cY=0$.
\end{proof}

We prove our main results:

\begin{Theorem}[\cref{71355712}]\label{74600988}
	Let $X$ be a divisorial noetherian scheme. Then the following conditions are equivalent:
	\begin{enumerate}
		\item\label{80221564} $\QCoh X$ satisfies Ab4*.
		\item\label{26021295} $\QCoh X$ has enough projectives.
		\item\label{49174160} $X$ is an affine scheme.
	\end{enumerate}
\end{Theorem}

\begin{proof}
	\cref{49174160}$\Rightarrow$\cref{26021295}: Since $\QCoh X\cong\Mod\Gamma(X,\cO_{X})$, it has enough projectives.
	
	\cref{26021295}$\Rightarrow$\cref{80221564}: See \cref{68007626}.
	
	\cref{80221564}$\Rightarrow$\cref{49174160}: By \cref{92619170} and \cref{50382042}, we have an equivalence $\QCoh X\isoto\Mod^{I}S$ that sends $\cO_{X}$ to an object isomorphic to $S$. Hence $\cO_{X}\in\QCoh X$ is a projective object, and we obtain
	\begin{equation*}
		H^{d}(X,-)\cong\Ext_{X}^{d}(\cO_{X},-)=0
	\end{equation*}
	for all integers $d\geq 1$. By Serre's criterion of affineness (\cite[Theorem~III.3.7]{MR0463157}), we conclude that $X$ is an affine scheme.
\end{proof}

\begin{Corollary}[\cref{66241253}]\label{44529212}
	Let $X$ be a scheme that contains a non-affine divisorial noetherian scheme as a closed subscheme. Then $\QCoh X$ does not satisfy Ab4*.
\end{Corollary}

\begin{proof}
	Let $Y\subset X$ be a closed subscheme with the stated property. It is shown in the proof of \cite[Corollary~3.9]{MR3361309} that the closed immersion $i\colon Y\into X$ induces the fully faithful functor $i_{*}\colon\QCoh Y\to\QCoh X$ whose essential image is the full subcategory $\cC$ of $\QCoh X$ consisting of all objects $\cN\in\QCoh X$ annihilated by the quasi-coherent subsheaf $\cI_{Y}\subset\cO_{X}$ corresponding to the closed subscheme $Y$. For every object $\cM\in\QCoh X$, the quotient object $\cM/\cM\cI_{Y}$ is largest among those belonging to $\cC$. Hence $\cC\subset\QCoh X$ is a closed subcategory by \cref{17471113}.
	
	By \cref{74600988}, $\QCoh Y$ does not satisfy Ab4*. Since condition Ab4* on a Grothendieck category is inherited by its closed subcategories, $\QCoh X$ does not satisfy Ab4*, either.
\end{proof}



\end{document}